\documentclass[11pt]{article}
\usepackage[margin=1.0in]{geometry}

\usepackage{amssymb,amsthm,amsmath,amsfonts,mathtools}

\usepackage{graphicx}
\usepackage{epstopdf}
\usepackage{subcaption}
\usepackage{adjustbox}
\usepackage{color}
\usepackage{tikz-qtree}
\usepackage{textcomp}
\usepackage[T1]{fontenc}

\usepackage{listings}
\usepackage[
    natbib=true,
    citestyle=alphabetic,
    style=alphabetic,
    maxnames=3,
    minnames=3,
    maxbibnames=99,
    backend=bibtex,
    urldate=iso8601,
    firstinits=true,
    isbn=false,
    doi=false,
    date=year,
    sorting=nty
]{biblatex}
\DeclareFieldFormat[article,inbook,incollection,inproceedings,patent,thesis,
  unpublished]{citetitle}{#1}
\DeclareFieldFormat[article,inbook,incollection,inproceedings,patent,thesis,
  unpublished]{title}{#1}
\renewbibmacro{in:}{%
  \ifentrytype{article}{}{\printtext{\bibstring{in}\space}}}
\renewcommand*{\labelalphaothers}{\textsuperscript{+}}

\usepackage[hidelinks]{hyperref}
\usepackage{breakurl}
\usepackage{url}

\newcommand{\rref}[1]{(\ref{#1})}
\newcommand{\secref}[1]{\S\ref{#1}}
\newcommand{\figref}[1]{figure~\ref{#1}}
\newcommand{\ie}{\textit{i.e.}}

\newcommand{\RR}{\mbox{\bf R}}

\newcommand{\reals}{{\mbox{\bf R}}}
\newcommand{\epi}{\mbox{\bf epi}}

\newcommand{\pf}{\lambda_{\text{pf}}}
\newcommand{\ee}[1]{e^{#1}}
\newcommand{\inv}[1]{#1^{-1}}
\newcommand{\Set}[2]{\{\,#1\,\mid\,#2\,\}}

\newcommand{\llcv}{log-log convex}
\newcommand{\llcc}{log-log concave}
\newcommand{\lla}{log-log affine}
\newcommand{\ullcv}{Log-log convex}

\makeatother

\begin{document}

\title{Disciplined Geometric Programming}

\author{
Akshay Agrawal \\ \texttt{
\small akshayka@cs.stanford.edu} \and
Steven Diamond \\ \texttt{\small diamond@cs.stanford.edu} \and
Stephen Boyd \\ \texttt{\small boyd@stanford.edu}
}
\date{December 12, 2018}

\maketitle

\begin{abstract}
We introduce {\llcv} programs, which are optimization problems with
positive variables that become convex when the variables, objective functions,
and constraint functions are replaced with their logs, which we refer to as a
log-log transformation.  This class of problems generalizes traditional
geometric programming and generalized geometric programming, and it includes
interesting problems involving nonnegative matrices.  We give examples of
{\llcv} functions, some well-known and some less so, and we develop an analog
of disciplined convex programming, which we call disciplined geometric
programming. Disciplined geometric programming is a subclass of {\llcv}
programming generated by a composition rule and a set of functions
with known curvature under the log-log transformation.  Finally, we describe an
implementation of disciplined geometric programming as a reduction in CVXPY
1.0.
\end{abstract}

\section{Introduction}\label{sec-intro}

\subsection{Geometric and generalized geometric programs}
A \textit{geometric program} (GP) is a nonlinear mathematical optimization
problem in which all the variables are positive and the objective
and constraint functions are either \textit{monomial
functions} or \textit{posynomial functions}. A monomial is any real-valued
function given by $x \mapsto cx_1^{a_1}x_2^{a_2} \cdots x_n^{a_n}$, where $x =
(x_1, x_2, \ldots, x_n)$ is a vector of positive real variables,
the coefficient $c$ is positive, and the exponents $a_i$ are real; 
a posynomial function is any sum of monomial functions.  A GP is
an optimization problem of the form
\begin{equation}\label{eqn-gp}
\begin{array}{ll}
\mbox{minimize} & f_0(x) \\
\mbox{subject to} & f_i(x) \leq 1, \quad i=1, \ldots, m\\
& g_i(x) = 1, \quad i=1, \ldots, p,
\end{array}
\end{equation}
where the functions $f_i$ are posynomials, the functions $g_i$ are monomials, and $x \in
\RR^n_{++}$ is the decision variable.
($\reals_{++}$ denotes the set of positive reals.)

The problem \rref{eqn-gp} is not convex, but it can be transformed to a convex
optimization problem by a well-known transformation.  We can make the change of
variables $u = \log x$ (meant elementwise) and take the logarithm of the
objective and constraint functions to obtain the equivalent problem
\begin{equation}\label{eqn-convex-gp}
\begin{array}{ll}
\mbox{minimize} & \log f_0(\ee{u}) \\
\mbox{subject to} & \log f_i(\ee{u}) \leq 0, \quad i=1, \ldots, m\\
& \log g_i(\ee{u}) = 0, \quad i=1, \ldots, p,
\end{array}
\end{equation}
which can be verified to be convex \cite[\S 4.5.3]{boyd2004}.
(The exponential $e^u$ is meant elementwise.)
Because GPs are reducible to convex programs, they can be solved efficiently
and reliably using any algorithm for convex optimization, such as
interior-point methods \citep{NN1994} or first-order methods
\cite{boyd2011distributed}. 
When all $f_i$ are monomials, the problem~(\ref{eqn-convex-gp}) reduces to 
a general linear program (LP), so GP is a generalization of LP.

Since its introduction four decades ago \citep{duffin1967},
geometric programming has found application in chemical engineering
\citep{clasen1984}, environment quality control \citep{greenberg1995}, digital
circuit design \citep{boyd2005}, analog and RF circuit design
\citep{hershenson2001opamp, li2004, xu2004}, transformer design \citep{jabr2005}, communication systems
\citep{kandukuri2002, chiang2005, chiang2007}, biotechnology \citep{marin2007,
vera2010}, epidemiology \citep{preciado2014}, optimal gas flow
\citep{misra2015optimal}, tree-water-network control \citep{sela2015}, and
aircraft design \citep{hoburg2014geometric, brown2018vehicle, saab2018robust}.
This list is far from exhaustive; for many other examples, see \S10.3 of
\citep{boyd2007}.

Evidently monomials and posynomials are closed under various operations.  For
example, monomials are closed under multiplication, division, and taking
powers, while posynomials are closed under addition, multiplication, and
division by monomials.  A \emph{generalized posynomial} is defined as a
function formed from monomials using the operations addition, multiplication,
positive power, and maximum.  Generalized posynomials, which include posynomials,
are also convex under a
logarithmic change of variable, after taking the log of the function.  It
follows that a \emph{generalized geometric program} (GGP), \ie, a problem of the
form \rref{eqn-gp}, with $f_i$ generalized posynomials and $g_i$ monomials,
transforms to a convex problem in~\rref{eqn-convex-gp} \citep[\S5]{boyd2007}, 
and therefore is tractable.

\subsection{{\ullcv} programs}
\begin{figure}
\centering
\includegraphics[width=0.30\textwidth]{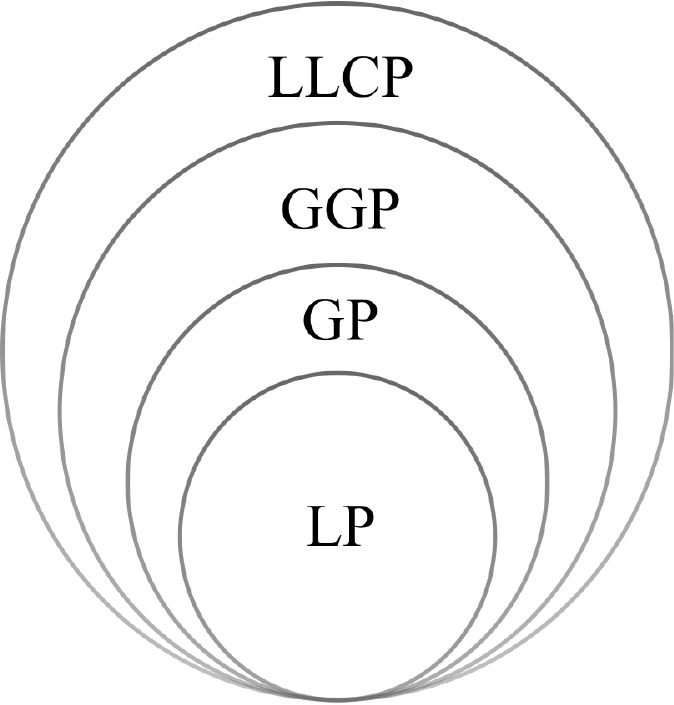}
\caption{Hierarchy of optimization problems.}
\label{fig-hierarchy}
\end{figure}

For a function $f: D \to \RR_{++}$, with $D \subseteq \RR^{n}_{++}$, we refer
to the function $F(u) = \log f(e^u)$, with domain
$\{u \mid e^u \in D \}$, as its \emph{log-log transformation}.
We refer to a
function $f$ as \emph{\llcv} if $F$ is convex, \emph{\llcc} if $F$ is concave,
and \emph{\lla} if $F$ is affine.  
As in convex analysis, we can consider the analog of extended-value extensions
\cite[\S 3.1.2]{boyd2004}:
we allow a {\llcv} function to take the value $+\infty$, and a {\llcc} function
to take the value zero,
which corresponds to $F$ taking the value $-\infty$.
A function is {\lla} if and only if it is a monomial;
posynomials and generalized posynomials are {\llcv}, but there are
{\llcv} functions that are not generalized posynomials (examples are given in
\secref{s-llcvx-simple-ex} and \secref{s-llcvx-matrix-ex}).

An optimization problem of the form \rref{eqn-gp}, with $f_i$ {\llcv} and
$g_i$ \lla, is called a \emph{{\llcv} program} (LLCP).  The set of
LLCPs is a strict superset of GGPs.
The hierarchy of LPs, GPs, GGPs, and LLCPs is shown in
figure~\ref{fig-hierarchy}.

{\ullcv}ity is also known as \textit{geometric convexity} or
\textit{multiplicative convexity}, since it is equivalent to convexity with
respect to the geometric mean (see \secref{s-prop}). \citet{montel1928} studied
the class of {\llcv} functions many decades ago, in the context of subharmonic
functions. More recently, \citet{niculescu2000} developed a theory of
inequalities derived from {\llcv}ity, parallel to the theory of convex functions,
\citet{forster2005spectral} studied the {\llcv}ity of certain operator polynomials,
and \citet{baricz2010} examined the log-log concavity of various univariate
probability distributions. See also \citep{doyle2006knee,
jarczyk2002mulholland, ozdemir2014note} for related work. 

Many functions can be well-approximated by {\llcv} functions \citep{boyd2007,
hoburg2016data, calafiore2018log}, but the lack of a coherent modeling
framework for LLCPs has hindered their use in practical applications. The point
of this paper is to close that gap.

\subsection{Domain-specific languages for convex optimization}
Disciplined convex programming (DCP) describes a subset of convex programs
generated by a single rule and a set of \emph{atoms}, functions with known
curvature (convex, concave, or affine) and monotonicity \citep{grant2006disciplined}. DCP is a natural
starting point for building a domain-specific language (DSL) for convex
optimization, \ie, a programming language that parses convex optimization
problems expressed in a human-readable form, rewrites them into
canonical forms, and supplies the lowered representations to numerical
solvers. By abstracting away solvers, DSLs make optimization accessible to
researchers and engineers who are not experts in the details of optimization
algorithms. Most DSLs for convex optimization have DCP as their foundation;
examples include CVX \citep{cvx}, CVXPY \citep{cvxpy, cvxpyrewriting},
Convex.jl \citep{convexjl}, and CVXR \citep{fu2017cvxr}.  For a survey of DSLs
for convex optimization, see \citep[\S1]{cvxpyrewriting}.
Some DSLs, like CVX and Yalmip \citep{yalmip}, can also
parse GPs and GGPs. There also exist DSLs specifically for GPs,
including GPKit \citep{gpkit} and GGPLAB \citep{mutapcic2006ggplab}. 
These software packages parse and rewrite GPs and GGPs.

In this paper, we introduce the analog of DCP for {\llcv}
problems.  We refer to our analog of DCP as \emph{disciplined geometric
programming} (DGP).  Like DCP, every disciplined geometric program is generated
by a single rule and a library of atoms. The class of disciplined
geometric problems is a subclass of {\llcv} problems (and of course depends on
the library of atoms), and, with a sensible atom library, a strict
superset of both geometric programming and generalized geometric programming.
In \secref{s-loglogcvx}, we characterize {\llcv}ity and
give many examples of {\llcv} functions, some obvious and some less
so; when appropriate, we also supply graph implementations
\citep{grant2008}. In \secref{s-dgp}, we
present DGP, along with a verification procedure that we articulate in terms of
mathematical expression trees.  We close in \secref{s-impl} by describing an
implementation of DGP as a reduction to disciplined convex programs in CVXPY
1.0.

\section{{\ullcv}ity}\label{s-loglogcvx}

\subsection{Properties}\label{s-prop}
\paragraph{Convexity with respect to the geometric mean.}
{\ullcv} functions obey a variant of Jensen's inequality: a function $f$ is
{\llcv} if and only if for all $x, y$ in the domain of $f$, and for each
$\theta \in [0, 1]$,
\[
f(x^\theta \circ y^{1 - \theta}) \leq f(x)^\theta f(y)^{1 - \theta},
\]
where $\circ$ is the Hadamard (elementwise) 
product and the powers are meant elementwise.

\begin{figure}
\centering
\adjustbox{valign=t}{\begin{minipage}[b]{0.25\textwidth}
\caption{Two {\llcv} functions and one {\llcc} function.}
\label{fig-log-log-plt}
\end{minipage}}%
\hfill
\adjustbox{valign=t}{\begin{minipage}[c]{0.73\textwidth}
\includegraphics[height=0.28\textheight]{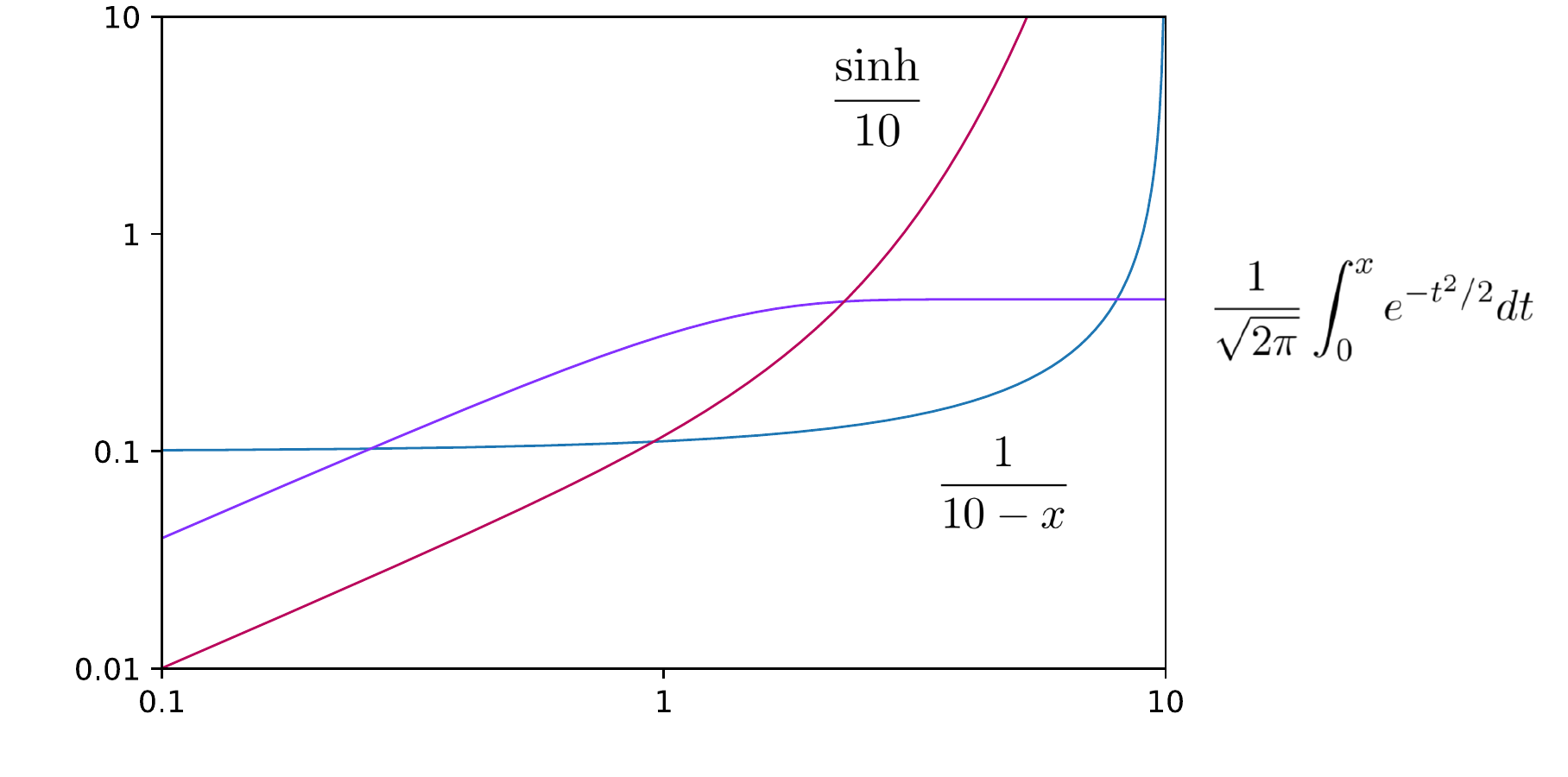}
\end{minipage}}
\end{figure}

\paragraph{Scalar {\llcv} functions.}
A scalar function $f : D \rightarrow \RR_{++}$, 
$D \subseteq \RR_{++}$, is {\llcv} if its graph
has positive curvature on a log-log plot, as shown in \figref{fig-log-log-plt}. If
$f$ is additionally twice-differentiable, then it is {\llcv} if and only if for
all $x\in D$,
\begin{equation}\label{eqn-second-order-scalar}
f''(x) + \frac{f'(x)}{x} \geq \frac{f'(x)^2}{f(x)}.
\end{equation}

\paragraph{Epigraph.}
If the set $\{u \mid e^u \in D\}$ is convex, $D \subseteq \RR_{++}^n$, we
say that $D$ is a \emph{log-convex} set. The domain of a {\llcv}
function $f$ is of course a log-convex set. Its epigraph
\[
\epi f = \Set{(x, t)}{f(x) \leq t}
\]
is also a log-convex set. The converse is true as well: if the epigraph
of a function is a log-convex set, then the function is {\llcv}.
These facts follow from the similar rules for convex functions and epigraphs
\citep[\S3.1.7]{boyd2004}.

\paragraph{Relationship to log-convexity.} {\ullcv} functions are related to
log-convex functions, which are real-valued functions $f$ for which $\log f$ is
convex \citep[\S3.5]{boyd2004}.  If $f$ is log-convex and nondecreasing in each of
its arguments, then its log-log transformation $F(u) = \log f (\ee{u})$ is
{\llcv}, as can be seen via the vector composition rule for convex functions
\citep[\S3.2.4]{boyd2004}.  Similarly, if $f$ is log-concave and nonincreasing in
its arguments, then its log-log transformation is \llcc. Since every positive
concave function is log-concave, it follows that every positive concave
function that is nonincreasing in its arguments is also \llcc.

In some cases, {\llcv}ity implies log-convexity. A function $f$ is log-convex
if and only if for all $x$ and $y$ in its domain and for each $\theta \in [0,
1]$,
\[
f(\theta x + (1 - \theta) y) \leq f(x)^ \theta f(y)^{1 - \theta}.
\]
In light of this fact and the AM-GM inequality, every nonincreasing {\llcv}
function is also log-convex, and every nondecreasing {\llcc} function is also
log-concave.

\paragraph{Partial minimization.} If $f$ is {\llcv} in the variables $x$
and $y$, and if $D$ is a log-convex set, then the function
\[
g(x) = \inf_{y \in D} f(x, y)
\]
is also {\llcv}. A similar result holds for {\llcc} functions: if $f(x, y)$
is {\llcc} and $D$ is a log-convex set, then $g(x) = \sup_{y \in D} f(x, y)$
is {\llcc}. These results are translations of identical results for convex
functions \citep[\S3.2.5]{boyd2004}.

\paragraph{Integration.} If $f : [0, a) \to [0, \infty)$ is continuous
and {\llcv} ({\llcc}) on $(0, a)$, then 
\[
x \mapsto \int_{0}^{x} f(t)dt
\]
is also {\llcv} ({\llcc}) on $(0, a)$ \citep{montel1928, niculescu2000}. As an
example, if $X$ is a real-valued random variable with a continuous {\llcc}
density $f$ defined on $[0, a)$, then the probability that $X$ lies between $0$
and some $x \in (0, a)$ is a {\llcc} function of $x$.  Several common
distributions, including the Gaussian, Gibrat, and the Student's
$t$, have {\llcc} densities \citep[\S5]{baricz2010}.

\subsection{Composition rule}\label{s-comp}
A basic result of convex analysis is that a nondecreasing
convex function of a convex function is convex.
(Similarly, a nonincreasing convex function of a concave function is
convex.) These results, along with similar ones for concave functions,
are special cases of just one result on the curvature
of function compositions, and it is on this single result that DCP is
based \citep[\S6.4]{grant2006disciplined}. An analogous composition rule
holds for {\llcv} functions, which we provide in full generality
below. Its proof is an elementary exercise in convex analysis.

Suppose $h : D \to \RR_{++} \cup \{\infty\}$, $D \subseteq \RR^{k}_{++}$, is {\llcv},
nondecreasing in its $i$th argument for each $i$ in an index set $I \subseteq
\{1, 2, \ldots, k\}$, and nonincreasing in the arguments indexed by $I^c$. For
$i = 1, 2, \ldots, k$, let $g_i : D_i \subseteq \RR^n_{++} \to \RR_{++}$.
Let $f : \bigcap D_i \to \RR_{++} \cup \{\infty\}$ be given by
\[
f(x) = h(g_1(x), g_2(x), \ldots, g_k(x)).
\]
If $g_i$ is {\llcv} for $i \in I$  and {\llcc} for $i \in I^c$,
then the function $f$ is {\llcv}.

A symmetric result holds when $h : D
\to \RR_{+}$, $D \subseteq \RR^k_{+}$, is \llcc: If
$g_i$ is {\llcc} for $i \in I$ and {\llcv} for $i \in I^c$, then
$f(x) = h(g_1(x), \ldots, g_k(x))$ is \llcc.

\subsection{Some simple examples}\label{s-llcvx-simple-ex}
We have already seen that monomials are {\lla} and that posynomials
and generalized posynomials
are {\llcv}. In this section we provide several other examples of
{\llcv} and {\llcc} functions.

\paragraph{Product.}
The product $f(x_1,x_2)=x_1x_2$ is \lla,
since $F(u) = \log(\ee{u_1}\ee{u_2}) = u_1 + u_2$ is affine.
(This is also clear since $f$ is a monomial.)
It follows that the product of {\lla} functions is {\lla},
and (since the product is monotone increasing) the product of
{\llcv} functions is \llcv, and the product of
{\llcc} functions is \llcc.

\paragraph{Ratio.}
The ratio $f(x_1,x_2)=x_1/x_2$ is {\lla} (since it is a monomial),
increasing in its first argument and decreasing in its second argument.
It follows that the ratio of a \llcv~and a \llcc~function is \llcv,
and that the ratio of \llcc~and a \llcv~function is \llcc.

\paragraph{Power.} For $a \in \RR$, the function given by $x^a$
is {\lla} in $x$, since $\log(\ee{ax}) = ax$.  It follows that a
power of a {\lla} function is \lla.
For $a \geq 0$,
the power of a {\llcv} function is \llcv,
and the power of a {\llcc} function is \llcc.
For $a < 0$,
the power of a {\llcv} function is \llcc,
and the power of a {\llcc} function is \llcv.

\paragraph{Sum.} The function $f(x_1, x_2) = x_1 + x_2$ is {\llcv} since $F(u) =
\log(\ee{u_1}+ \ee{u_2})$ is convex. It follows that the sum of {\llcv}
functions is \llcv. Log-log concavity is \emph{not} in general
preserved under addition.

\paragraph{Max and min.}  The function $f(x) = \max_i x_i $ is \llcv,
and the function $f(x)=\min_i x_i$ is \llcc.
Since both are nondecreasing, it follows that the max of \llcv~functions is 
\llcv, and the min of \llcc~functions is \llcc.

\paragraph{Sum largest.} For $x \in \RR^n_{++}$, the sum of the $r$ largest
elements in $x$ is \llcv, since it can be represented as
$\max \Set{x_{i_1} + x_{i_2} + \cdots + x_{i_r}}{i_1 < i_2 < \cdots < i_r}$,
which is the max of a finite number of {\llcv} functions.

\paragraph{One-minus.} The function $f(x) = 1 - x$ with domain $(0, 1)$
is \llcc, as can be seen by noting that $f$ is concave
and decreasing in $x$, or by the fact that the second derivative of its log-log
transformation is negative. It is also decreasing in $x$, so we
conclude that if $g$ is \llcv, $f(g(x)) = 1 - g(x)$ is {\llcc}
(with domain $\Set{x}{g(x) < 1}$).

\paragraph{Difference.}
The function $f(x) = x_1-x_2$, with domain $\{x>0 \mid x_1 - x_2 >0\}$,
is \llcc, increasing in its first argument and decreasing in its second.
It follows that the difference of a {\llcc} function and a log-log
convex function (with obvious domain) is \llcc.

\paragraph{Geometric mean.} The geometric mean $f(x) = \left(\prod_{i=1}^{n}
x_i\right)^{1/n}$ is {\lla}, \ie, a monomial. The geometric mean of
{\llcv} functions is \llcv, and likewise for \llcc~functions.

\paragraph{Harmonic mean.} The harmonic mean
$f(x) = n\inv{(1/x_1 + 1/x_2 + \cdots + 1 / x_n)}$ is {\llcc},
since it is the reciprocal of a {\llcv} function.

\paragraph{$\ell_p$-norm.} The $\ell_p$-norm
$\|x\|_p = (|x_1|^p + |x_2|^p + \cdots + |x_n|^p)^{1/p}$, $p \geq 1$, is
{\llcv} for $x \in \RR^n_{++}$, since $\|x\|_p$ with the absolute values
removed is a posynomial raised to $1/p$.

\paragraph{Exponential and logarithm.} The function $f$ given by $f(x) =
\ee{x}$ for $x > 0$ is \llcv, since $F(u) = \log f(\ee{u}) = \ee{u}$,
which is convex. Similarly, the logarithm function restricted to $(1, \infty)$
is \llcc.

\paragraph{Entropy.} The function $f(x)=-x\log x$ with domain $(0,1)$ is
\llcc, as can be seen via the composition rule.

\paragraph{Functions with positive Taylor expansions.}  Suppose $f:\reals \to \reals$
is given by a power series $f(x) = a_0 + a_1 x + a_2 x^2 + \cdots$,
with $a_i \geq 0$ and radius of convergence $R$.
We restrict $f$ to the domain $(0,R)$.
Then $f$ is \llcv.
This is readily shown by noting that the partial sums
are posynomials, so $f$ is the pointwise limit of \llcv~functions.
As examples, the functions $\sinh$ and $\cosh$ restricted to $(0,
\infty)$, $\tan$, $\sec$, and $\csc$ restricted to $(0, \pi / 2)$, $\arcsin$
restricted to $(0, 1]$, and $\log((1 + x)/(1 - x))$ restricted to $(0, 1)$ are
all \llcv.

\paragraph{Complementary CDF of a log-concave density.} The complementary
cumulative distribution function (CCDF) of a log-concave density is
{\llcc}. This follows from the fact that the CCDF of a log-concave density is
log-concave \citep[\S3.5.2]{boyd2004} and nonincreasing. As an example,
the CCDF of a Gaussian
\[
x \mapsto \frac{1}{\sqrt{2\pi}} \int_{x}^{\infty} e^{-t^2/2} dt
\]
is {\llcc} on $(0, \infty)$.
The densities of many common distributions, including the uniform,
exponential, chi-squared, and beta distributions, are log-concave. For
several other examples, see \citep[Table 1]{bagnoli2005log}.

\paragraph{Gamma function.} The Gamma function
\[
\Gamma(x)  = \int_{0}^{\infty}t^{x - 1}\ee{-t}dt
\]
is log-convex and nondecreasing for $x \geq 1$ \citep[\S3.5]{boyd2004}.
Hence, the restriction $\Gamma |_{[1, \infty)}$ is \llcv.

\subsection{Functions of positive matrices}\label{s-llcvx-matrix-ex}

In the following exposition, all inequalities should be interpreted
elementwise. For any two vectors $x, y$ in $\RR^n$, $x \leq y$ if and only
if the entries of $y - x$ are all nonnegative, and for any two matrices $A, B \in
\RR^{m \times n}$, we write $A \leq B$ to mean that the entries of $B - A$ are
nonnegative. Similarly, $x < y$ means that the entries of $y - x$ are positive,
and likewise for $A < B$. If $A > 0$, we will say that $A$ is a positive matrix.

Let
$\RR^{m \times n}_{++}$ denote the set of positive $m$-by-$n$
matrices. The log-log transformation of a function $f : D \subseteq \RR^{m \times n}_{++} \to \RR^{p
\times q}_{++}$ is $F(U) = \log f(\ee{U})$,
defined on $\Set{U}{e^U \in D}$, where the logarithm and exponential are meant
elementwise. We say that $f$ is {\llcv} if $F$ is convex with respect
to $\leq$, \ie, if for any $U$, $V$ in the domain of $F$, $\theta \in [0, 1]$
\[
F(\theta U + (1 - \theta) V) \leq \theta F(U)+ (1 - \theta) F(V).
\]
Equivalently, $f$ is {\llcv} if for any $X, Y \in D$, $\theta \in [0, 1]$,
\[
f(X^\theta \circ Y^{1 - \theta}) \leq f(X)^\theta f(Y)^{1 - \theta},
\]
where $\circ$ denotes the Hadamard product and the powers are meant elementwise.
Informally, we say that $f$ is {\llcv} if $f(X)$ has {\llcv}
entries for each $X \in D$.

Of course, the trace of a positive matrix and the product of positive matrices
are both {\llcv} functions. More interesting is the link between log-log
convexity and the Perron-Frobenius theorem, which states, among other
things, that every positive square matrix has a positive eigenvalue equal to its
spectral radius. We provide a few examples below.

\paragraph{Spectral radius.}
Let $X \in \mathbf{R}^{n \times n}$ have positive entries. The Perron-Frobenius
theorem states that $X$ has a positive real eigenvalue $\pf$ equal to its
spectral radius, \ie, the magnitude of its largest eigenvalue. It turns out
that $\pf$ is a {\llcv} function of $X$. This can be
seen by the fact that
\[
\pf = \min \Set{\lambda}{Xv \leq \lambda v \quad \textnormal{for
some } v > 0},
\]
where the inequalities are elementwise, which implies that $\pf \leq \lambda$
if and only if
\[
\sum_{j=1}^{n} X_{ij}v_j / \lambda v_i \leq 1, \quad i = 1, \ldots, n.
\]
The lefthand side of the above inequality is a posynomial in $X_{ij}$, $v_i$,
and $\lambda$, hence the epigraph of $\pf$ is log convex. This result
is described in more detail in \citep[\S4.5.4]{boyd2004}. For related material,
see \citep{kingman1961, nussbaum1986, forster2005spectral, weitan2015pf}.

\paragraph{Eye-minus-inverse.}
Let $D$ be the set of
positive matrices in $\RR^{n \times n}$ with spectral radius $\rho(X)$ less
than $1$. The function $f : D \rightarrow \RR^{n \times n}$ given by
\[
f(X) = (I - X)^{-1}
\]
is {\llcv} in $X$, \ie, $f(X)$ has {\llcv} entries.  The function $f$ is
well-defined: for any square matrix $X \in D$ , the power series $I +  X + X^2
+ \cdots$ converges to $(I - X)^{-1}$.  One intuitive way to see that $f$ is
\llcv~is to note that every partial sum $s_n(X) = \sum_{i=0}^{n} X^{i}$ with $n
\geq
1$ has posynomial entries, and therefore is \llcv. 
Because $s_n \rightarrow f$, we obtain that $f$ is \llcv.

We can also prove that the function $f$ is {\llcv} by studying its
epigraph. Let $X > 0$ and $T$ be matrices. Then
\begin{equation}\label{eqn-epi-lhs}
(I - X)^{-1} \leq T
\end{equation}
if and only if there exists a matrix $Y \geq 0$ such that
\begin{equation}\label{eqn-epi-rhs}
I \leq Y - YX, \quad Y \leq T.
\end{equation}
The equivalence between \rref{eqn-epi-lhs} and \rref{eqn-epi-rhs} shows that
the epigraph of $f$ is log convex: the set of matrices $X$, $Y$, and $T$
satisfying \rref{eqn-epi-rhs} is log convex, and the epigraph of $f$ is the
projection of this set onto its first and third (matrix) coordinates. It is
clear that \rref{eqn-epi-lhs} implies \rref{eqn-epi-rhs}, for if $X$ and $T$
satisfy \rref{eqn-epi-lhs}, then $X$, $Y=(I - X)^{-1}$, and $T$ satisfy
\rref{eqn-epi-rhs}. For the other direction, assume that the matrices $X > 0$,
$Y \geq 0$, and $T$ satisfy \rref{eqn-epi-rhs}. Let $\pf  = \rho(X) > 0$ be the
Perron-Frobenius eigenvalue of $X$ and let $v > 0$ be a corresponding right
eigenvector. Multiplying both sides of \rref{eqn-epi-rhs} by $v$, we obtain
that $v \leq (1 - \pf)Tv$. This necessitates that $\pf = \rho(X) < 1$, which
together with the fact that $X > 0$ implies that $(I - X)^{-1}$ exists and is
positive. Multiplying both sides of \rref{eqn-epi-rhs} by $(I - X)^{-1}$ yields
\rref{eqn-epi-lhs}.

\paragraph{Resolvent.}
For any square matrix $X$ and any scalar $s > 0$ such that $s$ is not an eigenvalue
of $X$, the matrix $(sI - X)^{-1}$ is called the \textit{resolvent} of $X$. The
function $(X, s) \mapsto (sI - X)^{-1}$ is {\llcv} in both $s$ and
$X$ whenever $X$ has positive entries and $\rho(X) < s$. This can be seen by
writing $(s I - X)^{-1}$ as $\inv{s}\inv{(I - X/s)}$.

\section{Disciplined geometric programming}\label{s-dgp}
While it is intractable to determine whether an arbitrary mathematical program is
{\llcv}, it is easy to check if a composition of
\textit{atoms} (functions with known log-log curvature and monotonicity)
satisfies the composition rule given in \secref{s-comp}. This fact motivates
\textit{disciplined geometric programming} (DGP), a methodology
for constructing {\llcv} programs from a set of atoms. A problem
constructed via disciplined geometric programming is called a
\textit{disciplined geometric program}. If a problem is a disciplined geometric
program, we colloquially say that the problem is DGP.

Like DCP \citep{grant2006disciplined}, DGP has
two key components: an atom library and a grammar for composing atoms.
Every function appearing in a disciplined geometric program must be either an
atom or a grammatical composition of atoms; a composition is grammatical if it satisfies
the rule from \secref{s-comp}.  Concretely, a disciplined geometric program is
an optimization problem of the form
\begin{equation}\label{eqn-dgp}
\begin{array}{ll}
\mbox{minimize} & f_0(x) \\
\mbox{subject to} & f_i(x) \leq \tilde{f}_i(x), \quad i=1, \ldots, m\\
& g_i(x) = \tilde{g}_i(x), \quad i=1, \ldots, p,
\end{array}
\end{equation}
where the functions $f_i$ are \llcv, the functions $\tilde{f}_i$
are \llcc, the functions $g_i$ and $\tilde{g}_i$ are \lla, $x \in
\RR^n_{++}$ is the decision variable, and all the functions are grammatical compositions
of atoms. (A problem where the objective is to
maximize a {\llcc} function and the constraints are as in
\rref{eqn-dgp} is also a disciplined geometric program.) Clearly, every disciplined geometric
program is an LLCP, but the converse is not true. This is not a
limitation in practice because atom libraries are extensible (\ie, the class
of DGP is parameterized by the atom library), and because invalid compositions
of atoms can often be appropriately re-expressed.

DGP offers an easy-to-understand prescription for constructing a large class of
{\llcv} problems. If the product, power, sum, and max functions are taken as
atoms, then DGP is equivalent to generalized geometric
programming. If other functions from \secref{s-llcvx-simple-ex} and
\secref{s-llcvx-matrix-ex} are also included, then the set of disciplined
geometric programs becomes a strict superset of the set of GGPs. As we shall
see in \secref{s-impl}, DGP is easily supported in
a DCP-based DSL for optimization.  For these reasons, it seems sensible to
suggest that DGP might replace GPs in the optimization modeling toolbelt.

\begin{figure}
\adjustbox{valign=t}{\begin{minipage}[c]{0.45\textwidth}
\includegraphics[height=0.28\textheight]{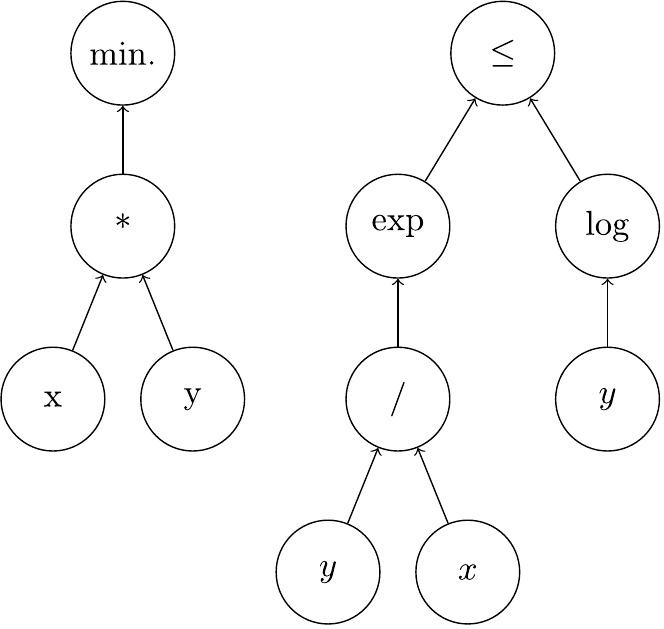}
\end{minipage}}%
\hfill
\adjustbox{valign=t}{\begin{minipage}[b]{0.55\textwidth}
\caption{Expression trees representing the optimization problem
(\ref{eqn-simple-ex}).
}
\label{fig-diff}
\end{minipage}}
\end{figure}

Verifying whether an optimization problem is DGP
involves representing the problem as a collection of mathematical
expression trees (one for the objective and one for each constraint), and
recursively verifying each expression tree. For example, the problem
\begin{equation}\label{eqn-simple-ex}
\begin{array}{ll}
\mbox{minimize} & xy \\
\mbox{subject to} & e^{y/x}  \leq  \log y
\end{array}
\end{equation}
can be represented by the expression trees shown in \figref{fig-diff};
assuming that the variables $x$ and $y$ are positive, this problem is an LLCP, but it
is neither a GP nor a GGP.

An expression tree for an objective is valid if its root is the minimize (maximize)
operator and the subtree rooted at its child is a valid {\llcv} ({\llcc}) composition of
atoms. A tree rooted at an atom is valid if the subtrees
rooted at its children are valid compositions of atoms, and if the composition
of the root with the subtrees of its children is grammatical.
Likewise, a tree for an inequality constraint is valid if
the left subtree is a valid {\llcv} composition of atoms, and the right subtree is
a valid {\llcc} composition of atoms. A tree for an equality constraint is
valid if both subtrees are {\lla} compositions of atoms. The recursion bottoms
out at the leaves of each tree, which are variables or constants. Leaves are
{\lla} provided that they are positive.

\section{Implementation}\label{s-impl}

We have implemented DGP in CVXPY 1.0, a Python-embedded, object-oriented DSL
for convex optimization \citep{cvxpyrewriting}. Our implementation, which is
available at
\begin{center}
  \url{https://www.cvxpy.org},
\end{center}
makes CVXPY 1.0 the first DSL for {\llcv} programming.

Our atom library includes a number of the functions
presented in \secref{s-llcvx-simple-ex} and \secref{s-llcvx-matrix-ex}, and
our implementation of DGP is a strict superset of generalized
geometric programming. CVXPY 1.0 can canonicalize any
DGP problem and furnish a solution to it, along with the optimal
dual values; it does this by reducing every DGP problem to a DCP problem,
canonicalizing and solving the DCP problem, and retrieving a solution to the
original problem.

\subsection{Canonicalization}
In CVXPY 1.0, canonicalization is facilitated by \texttt{Reduction} objects,
which rewrite problems of one form into equivalent problems of another form and
record how to retrieve a solution to the source problem from a solution to
the reduced-to problem. Canonicalizing DGP problems in CVXPY 1.0 is simple: we
first reduce each DGP problem to a DCP problem, after which we apply the DCP
canonicalization procedure.

We have added a class \texttt{Dgp2Dcp} that subclasses
\texttt{Reduction}. \texttt{Dgp2Dcp} accepts exactly those problems that
are DGP.  When applied to a problem, the \texttt{Dgp2Dcp} reduction recursively
replaces subexpressions with DCP log-log transformations or graph
implementations. For example, constants are replaced with their logarithms,
positive variables are replaced with unconstrained variables, products of two
expressions are replaced with sums of the log-log transformations of those
expressions, and sums of expressions are replaced with the
\texttt{log\_sum\_exp} of their canonicalized expressions. This procedure makes sense
because the log-log transformation of $f = h \circ g$ is equal to the
composition of the log-log transformations of $h$ and $g$.

Atoms like \texttt{eye\_minus\_inv} whose log-log transformations are
not DCP are replaced by their graph implementations (a graph implementation of
\texttt{eye\_minus\_inv} is given in \secref{s-llcvx-matrix-ex}). For
example, the expression \texttt{trace(eye\_minus\_inv(X))} would be
canonicalized to \texttt{trace(Y)}, together with the log-log transformation
of the constraint \texttt{Y U + I <= Y}, where \texttt{U} is a variable
representing \texttt{log X}.

\subsection{Solution retrieval}
When a DGP problem $\mathcal P_1$ is reduced to a DCP problem $\mathcal P_2$, for each variable
in $\mathcal P_1$, a variable representing its logarithm is instantiated
in $\mathcal P_2$.  Given a solution to $\mathcal P_2$, \ie, an assignment of numeric values to
variables, we recover a solution to $\mathcal P_1$ by exponentiating the values of the
variables in $\mathcal P_2$ and assigning the results to the corresponding variables in
$\mathcal P_1$. When $\mathcal P_2$ is unbounded, $\mathcal P_1$ is unbounded as well, in which case the
optimal value of the optimization problem is $0$ (if $\mathcal P_1$ is a minimization
problem) or  $+\infty$ (if $\mathcal P_1$ is a maximization problem).
Similarly, $\mathcal P_1$  is infeasible when $\mathcal P_2$ is infeasible.

The optimal dual values of $\mathcal P_1$ are the same as those of
$\mathcal P_2$. Under certain assumptions, the optimal dual values of $\mathcal P_2$ represent
fractional changes in the optimal objective given fractional changes in the
constraints \citep[\S3.3]{boyd2007}.

\subsection{Examples}
\paragraph{Hello, World.} Below is an example of how to use
CVXPY 1.0 to specify and solve the DGP problem (\ref{eqn-simple-ex}),
meant to highlight the syntax of our modeling language. A more interesting
example is subsequently presented.

\begin{lstlisting}[xleftmargin=.05\textwidth, xrightmargin=.2\textwidth]
import cvxpy as cp

x = cp.Variable(pos=True)
y = cp.Variable(pos=True)
objective_fn = x * y
objective = cp.Minimize(objective_fn)
constraints = [cp.exp(y/x) <= cp.log(y)]
problem = cp.Problem(objective, constraints)
problem.solve(gp=True)
print("Optimal value: ", problem.value)
print("x: ", x.value)
print("y: ", y.value)
print("Dual value: ", constraints[0].dual_value)
\end{lstlisting}
The optimization problem \texttt{problem} has two scalar variables, \texttt{x} and
\texttt{y}. For a problem to be DGP, every optimization variable must be
declared as positive, as done here with \texttt{pos=True}. The objective is to
minimize the product of $x$ and $y$, which is neither convex nor concave but is log-log affine,
since the product atom is log-log affine.  Every atom is an \texttt{Expression}
object, which may in turn have references to other \texttt{Expression}s; \ie,
each \texttt{Expression} represents a mathematical expression tree.
In line \texttt{7}, the \texttt{Expressions} are represented using
three atoms: ratio (\texttt{/}), \texttt{exp}, and
\texttt{log}. Also in line \texttt{7},
\texttt{exp(y/x)} is constrained to be no larger than \texttt{log(y)} via the
relational operator \texttt{<=},
which constructs a \texttt{Constraint} object linking two \texttt{Expressions}.
Line \texttt{8} constructs but does not solve 
\texttt{problem}, which encapsulates the expression trees for the objective and
constraints. The problem is DGP (which can be verified by asserting
\texttt{problem.is\_dgp()}), but it is not DCP (which can be verified by
asserting \texttt{not problem.is\_dcp()}).
Line \texttt{9} canonicalizes and solves \texttt{problem}. The optimal value of
the problem, the values of the variables, and the optimal dual value are
printed in lines \texttt{10}-\texttt{13}, yielding the following output.
\begin{lstlisting}[xleftmargin=.05\textwidth, xrightmargin=.2\textwidth,
numbers=none]
Optimal value:  48.81026898447343
x:  11.780089932635645
y:  4.143454698868564
Dual value:  2.843059917747706
\end{lstlisting}

As this code example makes clear, users do not need to know how
canonicalization works. All they need to know is how to construct DGP
problems. Calling the \texttt{solve} method on a \texttt{Problem} instance with
the keyword argument \texttt{gp=True} canonicalizes the
problem and retrieves a solution. If the user forgets to type
\texttt{gp=True} when her problem is DGP (and not DCP), a helpful error
message is raised to alert her of the omission.

\paragraph{Perron-Frobenius matrix completion.}
We have implemented several functions of positive matrices
as atoms, including the trace, product, sum, Perron-Frobenius eigenvalue, and
eye-minus-inverse.  As an example, we can use CVXPY 1.0 to formulate and solve
a Perron-Frobenius matrix completion problem. In this problem, we are given
\textit{some} entries of an elementwise positive matrix $A$,
and the goal is to choose the missing entries so as to minimize
the Perron-Frobenius eigenvalue or spectral radius.
Letting $\Omega$ denote the set of indices $(i, j)$ for which $A_{ij}$ is
known, the optimization problem is
\begin{equation}\label{eqn-pf-matrix-completion}
\begin{array}{ll}
\mbox{minimize} & \pf(X) \\
\mbox{subject to} & \prod_{(i, j) \not\in \Omega} X_{ij} = 1 \\
& X_{ij} = A_{ij}, \, (i, j) \in \Omega,
\end{array}
\end{equation}
which is an LLCP.  Below is an implementation of the problem
\rref{eqn-pf-matrix-completion}, with specific problem data
\begin{equation}\label{eqn-prob-data}
\centering
A = \begin{bmatrix}
1.0 & ? &  1.9 \\
? & 0.8 &  ? \\
3.2 & 5.9&  ?
\end{bmatrix},
\end{equation}
where the question marks denote the missing entries.

\begin{lstlisting}[xleftmargin=.05\textwidth, xrightmargin=.2\textwidth]
import cvxpy as cp

n = 3
known_value_indices = tuple(zip(*[[0, 0], [0, 2], [1, 1], [2, 0], [2, 1]]))
known_values = [1.0, 1.9, 0.8, 3.2, 5.9]
X = cp.Variable((n, n), pos=True)
objective_fn = cp.pf_eigenvalue(X)
constraints = [
  X[known_value_indices] == known_values,
  X[0, 1] * X[1, 0] * X[1, 2] * X[2, 2] == 1.0,
]
problem = cp.Problem(cp.Minimize(objective_fn), constraints)
problem.solve(gp=True)
print("Optimal value: ", problem.value)
print("X:\n", X.value)
\end{lstlisting}
Executing the above code prints the below output.

\begin{lstlisting}[xleftmargin=.05\textwidth, xrightmargin=.2\textwidth,
numbers=none]
Optimal value:  4.702374203221535
X:
[[1.         4.63616907 1.9       ]
 [0.49991744 0.8        0.37774148]
 [3.2        5.9        1.14221476]]
\end{lstlisting}

\clearpage
\printbibliography{}
\end{document}